\documentclass[11pt,reqno]{amsart}
\usepackage{amssymb,verbatim,amscd}
\usepackage{epsfig}
\usepackage[cp1251]{inputenc}\usepackage[russian]{babel}
\let\phi\varphi \let\le\leqslant \let\ge\geqslant

\let\opn\operatorname 

\newcommand{\CC }{{ \mathbb C {\,} }}\newcommand{\RR }{{ \mathbb R {\,} }}
\newcommand{\ZZ }{{ \mathbb Z {\,} }}
\newcommand{\QQ }{{ \mathbb Q {\,} }}
\newcommand{\CP }{{ \mathbb C P }}   

  \theoremstyle{plain}                          
  \newtheorem{THM}{Теорема}[section]            
  \newtheorem{LEM}[THM]{Лемма}
  \newtheorem{PROP}[THM]{Предложение}
  \newtheorem{COR}[THM]{Следствие}
  \newtheorem*{COR*}{Следствие}

  \theoremstyle{definition}                    
  \newtheorem{EXAM}[THM]{Пример}

  \theoremstyle{remark}                        
  \newtheorem{REM}[THM]{Замечание}

\numberwithin{equation}{section}

\let\ov\overline
\def\Nw{\operatorname{\mathit{Nw}}}

\def\nnnn{{\mathbf a}}

\def\SO{\mathrm{SO}}
\def\SU{\mathrm{SU}}
\def\Sp{\mathrm{Sp}}
\def\oU{\mathrm{U}}
\def\oS{\mathrm{S}}

\let\Cal\mathcal
\usepackage[mathcal]{euscript}

\title{}
\thanks{Поддержано РФФИ, грант 10-01-00041a.}
\keywords{Однородная метрика Эйнштейна, Homogeneous Einstein metric}
\address{Научно-исследовательский институт системных исследований РАН, Москва,
117218, Нахимовский проспект, 36, кор.~1}
\email{mmgraev@niisi.msk.ru}

\begin{document}



\hyphenation{ал-ге-б-ры ал-ге-б-ра ал-ге-б-рой ал-ге-бр ал-ге-б-рах ал-ге-б-ра-ми}
\hyphenation{ас-со-ци-а-ти-в-но-сть}
\hyphenation{ас-со-ци-и-ро-ван-но-го}
\hyphenation{аф-фин-ной аф-фин-ный аф-фин-ное аф-фин-ная аф-фин-но-му аф-фин-ном аф-фин-но-го аф-фин-ных}
\hyphenation{вся-ко-го вся-ких вся-кий вся-кой}
\hyphenation{вы-пол-ня-ет-ся}
\hyphenation{вы-пу-к-лым вы-пу-к-лой вы-пу-к-лый вы-пу-к-лая вы-пу-к-ло-го вы-пу-к-ло-му}
\hyphenation{ге-о-ме-т-рия ге-о-ме-т-рии ге-о-ме-т-р-ий ге-о-ме-т-р-и-ей ге-о-ме-т-рию}
\hyphenation{го-мо-те-ти-че-с-ких го-мо-те-ти-че-с-кие го-мо-те-ти-че-с-кий го-мо-те-ти-че-с-ко-го}
\hyphenation{ги-пер-пло-с-ко-сти ги-пер-пло-с-ко-сть ги-пер-пло-с-ко-ст-ей ги-пер-пло-с-ко-стях}
\hyphenation{гру-п-па гру-п-пы}
\hyphenation{дей-ст-ви-я-ми дей-ст-ви-ем дей-ст-вие дей-ст-вий дей-ст-вия дей-ст-вию}
\hyphenation{до-пус-ти-мос-ти}
\hyphenation{ев-к-ли-до-вых ев-к-ли-до-во ев-к-ли-дов ев-к-ли-до-ву ев-к-ли-до-ва}
\hyphenation{за-мы-ка-ния за-мы-ка-нию за-мы-ка-ние за-мы-ка-ни-ем}
\hyphenation{Зи-бен-та-ля}
\hyphenation{изо-мо-р-фи-з-мом изо-мо-р-фи-зм изо-мо-р-фи-з-ма изо-мо-р-фи-з-м-ов изо-мо-р-фи-з-м-а-ми изо-мо-р-фи-з-мам}
\hyphenation{изо-тро-пии}
\hyphenation{ин-ва-ри-а-нт-ных ин-ва-ри-ант-ные ин-ва-ри-ант-ный ин-ва-ри-ант-ная ин-ва-ри-ант-ное ин-ва-ри-ант-но-го ин-ва-ри-ант-ным ин-ва-ри-ант-ной}
\hyphenation{ин-ва-ри-ан-т-ны ин-ва-ри-ан-тен ин-ва-ри-ант-но}
\hyphenation{каж-до-го каж-до-му каж-д-ый каж-д-ая каж-д-ой каж-д-ых}
\hyphenation{ква-д-ра-тич-ная ква-д-ра-тич-ной ква-д-ра-тич-ную ква-д-ра-тич-ных ква-д-ра-тич-ным}
\hyphenation{ква-д-ра-т-ных ква-д-ра-т-ный}
\hyphenation{ком-па-к-т-но-го ком-па-к-т-но-му ком-па-к-т-ном ком-па-к-т-ным}
\hyphenation{ком-па-к-т-но ком-па-к-т-ен ком-па-к-т-на ком-па-к-т-ны}
\hyphenation{ком-па-к-ти-фи-ка-цию ком-па-к-ти-фи-ка-ция ком-па-к-ти-фи-ка-ции ком-па-к-ти-фи-ка-ци-ей}
\hyphenation{ком-п-лек-си-фи-ка-цию ком-п-лек-си-фи-ка-ция ком-п-лек-си-фи-ка-ции ком-п-лек-си-фи-ка-ци-ей ком-п-лек-си-фи-ка-ци-я-ми ком-п-лек-си-фи-ка-ций}
\hyphenation{ком-п-лек-с-ных ком-пле-кс-ны-ми}
\hyphenation{ком-по-нен-та ком-по-нен-те ком-по-нен-ты ком-по-нент ком-по-нен-той ком-по-нен-там ком-по-нен-тами}
\hyphenation{ко-неч-ным}
\hyphenation{клас-си-чес-ки-ми клас-си-чес-ким клас-си-чес-ких клас-си-чес-кий клас-си-чес-ко-го клас-си-чес-кой}
\hyphenation{ло-каль-но}
\hyphenation{ма-к-си-маль-но-го ма-к-си-маль-но-му ма-к-си-маль-но ма-к-си-маль-ны ма-к-си-маль-ный ма-к-си-маль-ные ма-к-си-маль-ных ма-к-си-маль-ным ма-к-си-маль-ном}
\hyphenation{ма-к-си-ма-лен}
\hyphenation{ме-т-ри-ки ме-т-ри-ка ме-т-ри-к-ой ме-т-рик ме-т-ри-ка-ми}
\hyphenation{мно-го-гран-ни-к мно-го-гран-ни-ка мно-го-гран-ни-ка-ми мно-го-гран-ни-кам мно-го-гран-ни-ком мно-го-гран-ни-ков мно-го-гран-ни-ку}
\hyphenation{мно-го-гран-ни-ки}
\hyphenation{мно-го-об-ра-зия мно-го-об-ра-зи-ях мно-го-об-ра-зие}
\hyphenation{мно-го-член мно-го-чле-ны мно-го-чле-на мно-го-чле-нов мно-го-чле-нам мно-го-чле-ном мно-го-чле-на-ми}
\hyphenation{мно-же-ст-в мно-же-ст-во мно-же-ст-ва мно-же-ст-ву мно-же-ст-ве мно-же-ст-вом мно-же-ст-вам мно-же-ст-ва-ми}
\hyphenation{Мин-ко-в-с-ко-го}
\hyphenation{на-и-мень-шее на-и-мень-шая на-и-мень-ший на-и-мень-шую на-и-мень-ше-го}
\hyphenation{не-ко-то-то-го не-ко-то-то-ой не-ко-то-то-рых}
\hyphenation{не-при-во-ди-мые не-при-во-ди-мых не-при-во-ди-мый не-при-во-ди-мо не-при-во-ди-мо-го}
\hyphenation{не-у-ни-мо-ду-ляр-но-го}
\hyphenation{од-но-кра-т-ным}
\hyphenation{од-но-ро-д-ных од-но-ро-д-ное од-но-ро-д-ном од-но-ро-д-но-го од-но-ро-д-ным од-но-ро-д-ные}
\hyphenation{од-но-свя-з-ной од-но-свя-з-ная од-но-свя-з-н-ых}
\hyphenation{опе-ра-то-р-ов опе-ра-то-ры опе-ра-то-ра опе-ра-то-ру опе-ра-то-рам опе-ра-то-ра-ми}
\hyphenation{оп-ре-де-ле-ние оп-ре-де-ле-ния оп-ре-де-ле-нию}
\hyphenation{от-но-си-тель-но}
\hyphenation{ото-б-ра-же-ний ото-б-ра-же-ни-я-ми ото-б-ра-же-ни-ем}
\hyphenation{пара-ме-три-за-ция пара-ме-три-за-ции пара-ме-три-за-ци-ей}
\hyphenation{пи-ра-ми-д пи-ра-ми-ды пи-ра-ми-да пи-ра-ми-де пи-ра-ми-дах пи-ра-ми-ду пи-ра-ми-дой пи-ра-ми-да-ми пи-ра-ми-дам}
\hyphenation{пло-с-кос-ти плос-кос-ть плос-кос-тей плос-кос-тями плос-кос-тью}
\hyphenation{пе-ре-ста-но-в-кой пе-ре-ста-но-в-ка пе-ре-ста-но-в-ки пе-ре-ста-но-вок пе-ре-ста-нов-кой пе-ре-ста-нов-ка-ми}
\hyphenation{по-д-мно-же-ст-ва-ми по-д-мно-же-ст-вo по-д-мно-же-ст-вoм по-д-мно-же-ст-вам}
\hyphenation{по-д-про-ст-ра-н-ст-во по-д-про-ст-ра-н-ст-в по-д-про-ст-ра-н-ст-ва по-д-про-ст-ра-н-ст-ва-ми по-д-про-ст-ра-н-ст-вам по-д-про-ст-ра-н-ст-вом}
\hyphenation{по-ло-жи-тель-но по-ло-жи-тель-но-го}
\hyphenation{по-лу-про-ст-ой по-лу-про-ст-ую по-лу-про-ст-ые по-лу-прост по-лу-про-ста по-лу-про-сты -ых по-лу-про-ст-о-го по-лу-про-ст-ой}
\hyphenation{по-ст-ро-ить}
\hyphenation{по-с-ре-д-ст-вом}
\hyphenation{пре-д-ста-в-ле-ния пре-д-ста-в-ле-ние пре-д-ста-в-ле-ний пре-д-ста-в-ле-нием пре-д-ста-в-ле-ни-ям пре-д-ста-в-ле-ни-я-ми}
\hyphenation{пре-дло-же-ни-я пре-дло-же-ни-е пре-дло-же-ни-й}
\hyphenation{про-ст-р-а-н-ств про-ст-р-а-н-ст-вом про-ст-р-а-н-ст-во про-ст-р-а-н-ст-ву про-ст-р-а-н-ст-ве про-ст-р-а-н-ст-ва про-ст-р-а-н-ст-вах про-ст-р-а-н-ст-ва-ми}
\hyphenation{про-хо-дя-щие про-хо-дя-щий про-хо-дя-щих про-хо-дя-щим про-хо-дя-щи-ми про-хо-дя-ще-го про-хо-дя-щ-ую про-хо-дя-щ-ая}
\hyphenation{про-ти-во-по-лож-ных про-ти-во-по-лож-ный про-ти-во-по-лож-ные про-ти-во-по-лож-ная про-ти-во-по-лож-ное про-ти-во-по-лож-ной про-ти-во-по-лож-но-го про-ти-во-по-лож-ным про-ти-во-по-лож-ны-ми}
\hyphenation{ра-бо-те ра-бо-та ра-бо-ты}
\hyphenation{рас-ши-ре-ни-ем}
\hyphenation{ра-ци-о-наль-ной ра-ци-о-наль-ный ра-ци-о-наль-ное ра-ци-о-наль-но-го ра-ци-о-наль-ных ра-ци-о-наль-ным ра-ци-о-наль-ны-ми}
\hyphenation{ре-зуль-та-те ре-зуль-тат ре-зуль-та-тов ре-зуль-та-ту ре-зуль-та-том ре-зуль-та-та-ми}
\hyphenation{ри-ма-но-ву ри-ма-но-ва ри-ма-но-во ри-ма-но-вых ри-ма-нов}
\hyphenation{рич-чи}
\hyphenation{сим-ме-т-ри-чес-ких сим-ме-т-ри-чес-кие сим-ме-т-ри-чес-кое сим-ме-т-ри-че-с-кая сим-ме-т-ри-че-с-ки-ми}
\hyphenation{сим-пли-ци-а-ль-ный сим-пли-ци-а-ль-ная сим-пли-ци-а-ль-но-го сим-пли-ци-а-ль-ным сим-пли-ци-а-ль-ных}
\hyphenation{си-сте-ма си-сте-ме си-сте-мой си-сте-мы си-стем}
\hyphenation{сле-до-ва-тель-но}
\hyphenation{со-от-ве-тст-в-ен-но}
\hyphenation{со-от-ве-тст-в-у-ю-щ-их}
\hyphenation{со-по-ста-вить со-по-ста-вив со-по-ста-вим}
\hyphenation{со-по-ста-в-лять со-по-ста-в-ля-ет-ся со-по-ста-в-ле-но со-по-ста-в-ле-ны со-по-ста-в-лен со-по-ста-в-лен-ный со-по-ста-в-лен-но-му со-по-ста-в-ля-е-мо-му со-по-ста-в-ля-е-мым}
\hyphenation{сжа-т-ых}
\hyphenation{ска-ляр-ную ска-ляр-ная ска-ляр-но-го ска-ляр-ных ска-ляр-ным ска-ляр-ны-ми ска-ляр-ные ска-ляр-ное ска-ляр-ный}
\hyphenation{спе-к-т-ром}
\hyphenation{то-ри-че-с-ко-го то-ри-че-с-кое то-ри-че-с-к-их}
\hyphenation{то-ч-но-стью то-ч-но-сть}
\hyphenation{тран-с-вер-с-а-ль-но}
\hyphenation{три-ан-гу-ля-ция три-ан-гу-ля-ции три-ан-гу-ля-ций три-ан-гу-ля-ци-ей три-ан-гу-ля-ци-я-ми}
\hyphenation{удо-в-ле-тво-рять удо-в-ле-тво-ря-ют удо-в-ле-тво-ря-ю-щие удо-в-ле-тво-ря-ю-щий удо-в-ле-тво-ря-ю-щих удо-в-ле-тво-ря-ю-щей удо-в-ле-тво-ря-ю-ще-го удо-в-ле-тво-ря-ю-щее}
\hyphenation{урав-не-ния урав-не-ние урав-не-ний урав-не-ни-ем урав-не-ни-я-ми}
\hyphenation{ус-ло-вие ус-ло-вия ус-ло-ви-ям ус-ло-ви-ем ус-ло-ви-я-ми ус-ло-ви-ям ус-ло-вию ус-ло-вий}
\hyphenation{ут-вер-ж-да-ет-ся ут-вер-ж-да-ет ут-вер-ж-да-ют}
\hyphenation{хо-ро-шей}
\hyphenation{фа-к-то-р-про-ст-ра-н-ст-ву фа-к-то-р-про-ст-ра-н-ст-во фа-к-то-р-про-ст-ра-н-ст-ва фа-к-то-р-про-ст-ра-н-ст-вом фа-к-то-р-про-ст-ра-н-ст-вам фа-к-то-р-про-ст-ра-н-ст-ва-ми фа-к-то-р-про-ст-ра-н-ств}
\hyphenation{фи-к-си-ро-ван-ной фи-к-си-ро-ван-ный фи-к-си-ро-ван-ная фи-к-си-ро-ван-ную}
\hyphenation{фи-к-си-ру-ем}
\hyphenation{филь-тра-ция}
\hyphenation{фла-го-вых фла-го-вое фла-го-во-го фла-го-вым фла-го-вом}
\hyphenation{функ-ция функ-ции функ-ций функ-ци-ей функ-ци-я-ми}
\hyphenation{Эйн-ш-тей-на}
\hyphenation{эк-ви-ва-лен-т-но эк-ви-ва-лен-т-на эк-ви-ва-лен-т-ны эк-ви-ва-лен-тен}
\hyphenation{эк-ви-ва-лен-т-ное эк-ви-ва-лен-т-ные эк-ви-ва-лен-т-ную эк-ви-ва-лен-т-ной эк-ви-ва-лен-т-но-го эк-ви-ва-лен-т-ных эк-ви-ва-лен-т-ным эк-ви-ва-лен-т-ны-ми}
\hyphenation{яв-ля-е-тся}





\noindent 18 oct. 2011 7.07

 \begin{center}{}\bf
On invariant Einstein metrics on K\"ahler homo\-geneous spaces
\\
$SU_4/T^3$,
$G_2/T^2$,
$E_6/T^2\cdot (A_2)^2$, $E_7/T^2\cdot A_5$, $E_8/T^2\cdot E_6$, $F_4/T^2\cdot A_2$
\\[1ex]\rm
Michail M. Graev
\end{center}

\noindent
{\footnotesize
{\sc Abstract.}
 \let\epsilon\varepsilon
 We study invariant Einstein metrics on the indicated homogeneous manifolds $M$, the corresponding
algebraic Einstein equations $E$, the associated with $M$ and $E$
Newton polytopes $P(M)$, and the integer volumes $\nu = \nu(P(M))$
of  it  (the Newton numbers). We show that $\nu = 80, 152, \dots,152$ respectively. 
It is claimed that the numbers $\epsilon = \epsilon(M)$
of complex solutions of $E$ equals $ \nu - 18, \nu - 18, \nu, \dots, \nu $.
The results are consistent with classification of non K\"ahler
invariant Einstein metrics on $G_2/T^2$ obtained recently by
Y.Sakane, A. Arvanitoyeorgos, and I. Chrysikos. We present also a
short description of all invariant complex Einstein metrics on 
$ SU_4/T^3 $. We prove existence of Riemannian non K\"ahler invariant
Einstein metrics on $G_2/T^2$-like K\"ahler homogeneous spaces 
$ E_6/T^2\cdot(A_2)^2 $, $ E_7/T^2\cdot A_5 $,  $E_8/T^2\cdot E_6 $,
$ F_4/T^2\cdot A_2 $, where $ T^2\cdot A_5 \subset  A_2\cdot A_5\subset
E_7 $ and some other results.

}

{\maketitle}
\vskip-1cm
\vskip-1cm

 \begin{center}{}\bf
Об инвариантных эйнштейновых метриках на кэлеровых
однородных пространствах
$SU_4/T^3$,
$G_2/T^2$,
$E_6/T^2(A_2)^2$, $E_7/T^2A_5$, $E_8/T^2E_6$, $F_4/T^2A_2$
\end{center}

\vskip-0.5cm
\vskip-0.5cm
\vskip-1cm
\renewcommand{\contentsname}{}
\setcounter{tocdepth}{2}
{\footnotesize\tableofcontents}
\vskip-1cm
\vskip-1cm




\section{Введение}

Инвариантные положительно определенные метрики Эйнштейна
на 12-мерном однородном пространстве $SU_4/T^3$
фактически были  классифицированы
\footnote{
Эта работа Ю.Сакане почти не содержит выделенных формулировок,
и в ней многое дается расположением материала.
Гиббонс, Лу и Поуп \cite{GLP}
оспаривают полноту ее изложения
в сходном вопросе о классификации для $Sp(2)/T^2$.
}
в статье Ю.Сакане \cite{Sa}.

\begin{THM}[\cite{Sa}]{}\label{THM:1}
На 12-мерном кэлеровом однородном пространстве $SU_4/T^3$
с точностью до изометрии и умножения на скаляр существуют три и только три
инвариантные метрики Эйнштейна $g$, не допускающие никакой инвариантной кэлеровой структуры.
Каждая из них имеет в группе Вейля, т.е. симметрической группе $S_4$,
нетривиальный стабилизатор, содержащий инволюцию $\sigma $ цикленного типа
$(12)$. Все метрики $g$ знакоопределены.
\end{THM}

Вместе с метрикой Кэлера--Эйнштейна получаются четыре попарно
неэк\-ви\-ва\-лен\-т\-ные
инвариантные
эйнштейновы метрики на $SU_4/T^3$, а вовсе не три, как утверждается в
\cite[таблица 1]{11}
со ссылкой на \cite{Sa}.
Заметим, что стабилизатор инвариантной кэлеровой метрики на $SU_4/T^3$
в группе $S_4$ меняет местами две противоположные камеры Вейля
(т.е. порожден инволюцией типа $(14)(23)$)
и тогда не содержит отражения $\sigma $.

Теорема~\ref{THM:1}
подсказала автору следующий результат,
анон\-си\-рованный в
\cite[при\-ме\-чание при корректуре]{2007}
и рассмотренный в \S8 диссертации \cite{2008}.

\begin{THM}{}\label{THM:2}
На $G_2/T^2$-подобных кэлеровых однородных пространствах
\begin{equation}{}\label{eq:1}
G_2/T^2,        \quad
E_6/T^2\cdot(A_2)^2, \quad
E_7/T^2\cdot A_5,     \quad
E_8/T^2\cdot E_6,     \quad
F_4/T^2\cdot A_2
\end{equation}
(где $T^2\cdot A_5 \subset A_2\cdot A_5  \subset E_7$)
существуют инвариантные положительно определенные метрики Эйнштейна $g$,
допускающие изометрию, соответствующую отражению
в группе Вейля системы корней типа $G_2$
(т.е группе диэдра порядка $12$)
и вследствие этого не допускающие никакой инвариантной кэлеровой структуры.
\end{THM}

А.Арванитойеоргос, И.Хрисикос и Ю.Сакане \cite{SACh}
доказали эту теорему независимо для случая
$12$-мерного пространства $G_2/T^2$.
Они также явно нашли инвариантные метрики Эйнштейна.
Из их формул следует:

\begin{THM}[\cite{SACh}]{}
На 12-мерном кэлеровом однородном пространстве $G_2/T^2$
с точностью до изометрии и умножения на скаляр существует
$8$ и только $8$ инвариантных метрик Эйнштейна, а именно:
\begin{enumerate}
\item[---]
две положительно определенные метрики, 
не допускающие никакой инвариантной кэлеровой структуры;
каждая из них лежит на $3$-элементной орбите группы Вейля
и удовлетворяет теореме~\ref{THM:2};

\item[---]
пять попарно неизометричных индефинитных метрик;
каждая из них  лежит на $6$-элементной орбите группы Вейля;

\item[---]
хорошо известная инвариантная метрика Кэлера-Эйнштейна
(единственная с точностью до изометрии).
\end{enumerate}
С точностью до умножения на скаляр на $G_2/T^2$
существует
$12$ знакоопределенных и $30$ индефинитных
инвариантных метрик Эйнштейна.
\end{THM}

Приближенные значения для положительно определенных метрик
первого вида,
найденные в \cite{SACh}
совпадают с вычисленными в \cite{2008}.

В \cite{SACh} на с. 18 выписаны приближенные значения только
для $14$ индефинитных ре\-ше\-ний, где $9,10,11,12$ и $13$-е
принадлежат пяти различным орбитам группы Вейля
(заметим, что центр этой $12$-членной группы состоит из двух элементов,
сохраняющих каждую инвариантную метрику). Это дает
$30=5\cdot 6$ индефинитных метрик с точностью до умножения на скаляр.
Других вещественных решений, порождающих $6$-эле\-мен\-т\-ные орбиты,
не существует, ибо (как показывает простое рассуждение)
при\-ве\-ден\-ное в \cite{SACh} на трех страницах алгеб\-ра\-и\-ческое уравнение
$84$-й степени для одной переменной имеет $30$ вещественных корней.

%
%

Перейдем к комплексным решениям уравнения Эйнштейна для
инвариантных метрик, рассматриваемых с точностью до умножения на комплексное
число.

По формуле (8.4) из \cite{2008} на $G_2/T^2$ должно быть $152-18$
комплексных решений с учетом алгебраических кратностей.
За вычетом $6$ решений Кэлера-Эйнштейна и еще $44$,
допускающих отражение из группы Вейля, остается
$$
84 = 152 - 18 - 50
$$
комплексных решений.
Значение $84$ было получено с помощью техники многогранника Ньютона.
без вычисления всех решений.
Теперь, что особенно приятно, оно подтверждено непосредственно
\cite{SACh}.
Вдобавок программа MAPLE сообщает,
что дискриминант трехстраничного уравнения $84$-й степени из
\cite{SACh} отличен от $0$.
Это подтверждает гипотезу \cite{2008} об однократности всех комплексных решений.

Читатель может обратиться к \cite{2006,2007} за определениями
многогранника Ньютона $\Nw(G/H)$ и целого числа Ньютона $\nu(G/H)$,
сопоставляемых связному риманову
одно\-родному пространству $G/H$ с компактной группой изотропии $H$
с простым спектром представления изотропии.
Как там, через $\varepsilon (G/H)$
обозначается число изолированных
 {\bf  комплексных}
решений алгебраического уравнения Эйнштейна для инвариантных метрик на $G/H$
(рассматриваемых с точностью до комплексного множителя)
с учетом алгебраических кратностей
этих решений.
{\bf Дефектом}
называется разность
$$
\delta (G/H) := \nu (G/H)- \varepsilon (G/H);
\qquad \delta (G/H)\ge0.
$$
Числа $\delta (G/H)$ и $\nu (G/H)$ можно находить и без явного
отыскания решений. А именно, $\nu$ является приведенным объемом
многогранника Ньютона, а $\delta $ иногда удается получать,
исследуя асимптотические решения, соответствующие
граням этого многогранника.
При $\nu>0$
асимтотические решения отсутствуют если и только если $\delta =0$;
если они отсутствуют или только изолированы,
то все настоящие решения изолированы.

Например, для $SU_4/T^3$ и $G_2/T^2$
можно явно найти все асимптотические решения
и убедиться, что они изолированы (ср. \cite{2007,2008}). 
Исследованием кратностей только асимптотических решений
мною доказано, что в этих случаях $\delta = 18$.

Как доказано в \cite[\S7]{2007}, для однородного пространства
$M=SU_4/T^3$ выполняется
$$
\nu = 80, \quad
\delta = 18.
$$
При этом на $M$ существует одно четырехкратное решение (это стандартная метрика Вана--Циллера)
и число попарно различных комплексных решений (здесь и далее с точностью до умножения на скаляр)
равно $\varepsilon - 3 = 80-18-3 = 59$.

Согласно \cite[\S8]{2008}, для пяти $G_2/T^2$-подобных
однородных пространств $M$ вида \eqref{eq:1} выполняется
$$
\nu = 152, \quad  \delta = \begin{cases}
18, & M= G_2/T^2,\\ 0, &\mbox{в остальных случаях}.
\end{cases}
$$
Во всей пяти случаях
число инвариантных метрик Кэлера-Эйнштейна равно половине от числа камер Вейля
системы $G_2$, т.е. $6$, а число комплексных решений, допускающих хотя бы
одно отражение из группы Вейля, равно $44$
(сверх того, все решения изолиро\-ва\-ны). 
В случае $M=G_2/T^2$ подсчет, основанный на предположении об однократности
остальных комплексных решений, и некоторые элементарные аргументы, дают
$k=14$ шестиэлементных орбит группы Вейля:
$$
k=\frac16(\nu-50- \delta )= \frac16( 152-50-18) = 14
$$

В заключение приведем теорему о дефекте из \cite{2006,2007}:

\begin{THM}
Пусть $G$ --- компактная простая группа Ли ранга $n\ge2$,
$T^{\,n}$ --- ее максимальный тор, $M=G/T^{\,n}$.
Тогда или $\nu(M)\ne \mathcal E(M)$, или $M=\SU_3/T^{\,2}$.
\end{THM}

Доказательство теоремы, данное в \cite{2007},
требует отдельного рассмотрения случаев
$SU_4/T^3$, $Sp_2/T^2$ (\cite[\S7]{2007}) и $G_2/T^2$
(см. \cite[\S8]{2008} или ниже, \S3).
Неравенство $\delta (G_2/T^2)>0$ можно доказать очень просто,
не вычисляя явно $\nu(G_2/T^2)$ и $\varepsilon (G_2/T^2)$,
а тем более, не решая уравнения Эйнштейна. 

Содержание этой статьи следующее. В \S2 приводится классификация
комплексных решений уравнения Эйнштейна для инвариантных метрик на
$SU_4/T^3$,
включающая классификацию Сакане вещественных метрик Эйнштейна.
В \S3 перепечатан раздел из диссертации \cite{2008},
касающийся пяти однородных пространств $M$ вида \eqref{eq:1}
и факторпространств $M/ \sigma $,
где $\sigma $ соответствует отражению из группы Вейля $W(G_2)$.
Вначале приводятся формулировки.
В \S\S 3.1--3.3 доказано, что $\delta (M/\sigma )=0$, $\nu (M/\sigma )=16$.
Для сравнения анализируются
фактор $SU_4/T^3$ по инволюции и
другие однородные пространства с общим многогранником Ньютона $P_{12}$
(трехмерная призма).
Для них $\delta  =0$, $\nu = 12$.
В \S 3.4 доказано, что $\nu(M)=152$, а в \S 3.6 -- что $\delta (G_2/T^2)>0$.
В добавлениях к \S 3 строятся 16 (комплексных)
эйнштейновых метрик на каждом пяти пространств $M/ \sigma $,
среди которых есть вещественные положительно определенные;
в случае $M=G_2/T^2$ получаются положительно определенные
метрики Сакане--Арванитойеоргоса--Хрисикоса.



\section{Классификация инвариантных эйнштейновых метрик на $SU_4/T^3$}
\label{sect:prim-pri-korr}

Далее под гомотетией понимается умножение метрики на скаляр.
Используются обозначения для метрик из \cite{2006,2007}.
Касательное расслоение к $SU_4/T^3$ расщепляется
на $6$ неприводимых инвариантных подрасслоений.
Отсюда каждая инвариантная метрика на $SU_4/T^3$ задается естественными
координатами $t_{ij}=t_{ji}$, $1\le i<j\le 4$.\quad
\footnote{
Содержание этого раздела в основном соответствует примечанию при корректуре
из \cite{2008}
}

Применяя полученную в \cite[\S7.3]{2007} формулу дефекта к пространству $M=\SU_4/T^3$,
на\-хо\-дим $\delta _M = 18$, откуда ${\mathcal E}(M) = 80-18 = 62$.
Как отмечалось, 
без учета кратностей
существует не $62$, а $59$
го\-мо\-те\-ти\-чес\-ких
классов (г.к.) инвариантных комплексных
метрик Эйнштейна в $M$. А именно, г.к.
двадцати девяти положительно определенных и тридцати мнимых метрик

 С точностью до гомотетии и изометрии на $M=\SU_4/T^3$
существуют
че\-ты\-ре
положительно определенные инвариантные метрики Эйнштейна $g_t^M$, а вовсе не
три,
как утверждается в \cite[таблица~1]{11}  со ссылкой на работу Ю.Сакане.
Достаточно их клас\-си\-фи\-ци\-ро\-вать с точностью до гомотетии и до
правого действия в $M$ группы Вейля $S_4$
(так как по теореме Онищика \cite{12} связная компонента группы
всех изометрий ин\-ва\-ри\-ант\-ной метрики на $M$
совпадает с
$\mathrm P\SU_4$).
Следующие четыре римановы метрики попарно неэквивалентны:

\par
(i) метрика Кэлера--Эйнштейна с 
координатами $t_{ab}=b-a$, $1\le a{<}b \le 4$;

\par
(ii) стандартная метрика Вана--Циллера 
с координатами $t_{ab} = 1$, $1\le a{<}b \le 4$;

\par
(iii) метрика Арванитойеоргоса: 
$t_{ij}=5$, $t_{k4}=3$ при ${\{i,j,k \}=\{1,2,3\}}$;

\par
(iv) риманова метрика вида
$t_{12} = {(2 + 3\cdot \theta)\,(7 - 12\cdot \theta)}/{(17 - 6\cdot \theta)}$,
$t_{34} = 2/3  + \theta$,
$t_{23}=t_{13}=1+ \eta $,
$t_{14}=t_{24}=1- \eta $,
где
$
\eta = (\theta + 2/3) \,\sqrt{{(5 + 3 \cdot \theta )}/{(34 - 12\cdot \theta)}}
$,
а $ \theta  = \theta _1 \approx - 0.028 $ --- действительный корень
приводимого над полем $\mathbb Q(l)$, 
где $$l^3 = 3\sqrt{57}\pm1 = 3\sqrt{A(6,2)}\pm1,$$
кубического уравнения
$
f(\theta )  = {2}/{27} + {8}/{3} \cdot\theta  + \theta ^3 = 0.
$

\par
Следующие три мнимые комплексные метрики попарно неэквивалентны:
метрики (v) и~(vi) вида (iv), отвечающие
мнимым корням $\theta = \theta _2$, $\theta _3$ многочлена $f(\theta )$,
где $|\theta _2/\theta_1| \approx 59$;
метрика~(vii)
с координатами
$t_{12}= \overline{t_{34}} = (4+\sqrt{-2})/3$ и
$t_{ij}=t_{k4} =1$ в остальных случаях.

\par
Из предложения 7.1 в \cite{2007} следует, что
каждую инвариантную комплексную метрику Эйн\-штей\-на в $M$
мож\-но перевести в одну и только одну из метрик (i)-(vii) действием
группы $S_4$ и ум\-но\-жением на скаляр,
т.е. существует
$59 = 29 + 30 = (12+1+4+12) + (12 + 12 +6)$
г.к.,
состоящих из метрик Эйнштейна.
Уравнение Эйнштейна 
\cite[Eq. (1.1)]{2007}
для метрик (i)-(vii)
можно проверить элементарно,
причем г.к. метрики (ii)
является его четырехкратным ре\-ше\-ни\-ем.
Осталось воспользоваться неравенством ${\mathcal E}(M) \le 62 = 59+3$.

\par
Выяснилось, что в упомянутой работе Ю.Сакане 
(см. \cite{Sa})
в явной но громоздкой форме построена метрика (iv)
и указано, что этот случай ранее рассматривал R.Senda.



\par
Судя по \cite[таблица~1]{11}, существование инвариантных
''римановых некэлеровых''
метрик Эйнштейна в $G_2/T^2$ было открытым вопросом.
Теперь
они найдены
на всех флаговых пространствах с си\-сте\-мой $T$-корней типа $G_2$,
т.е.
на $G_2/T^2$,
$E_6/T^2(A_2)^2$, $E_7/T^2A_5$, $E_8/T^2E_6$, $F_4/T^2A_2$.
Следующий текст продолжает \cite{2007} и использует те же обозначения.



\section{Флаговые пространства с системой $T$-корней типа $G_2$}
\label{sect:8}

Пусть $M$ --- флаговое пространство $SU_4/T^3$ или $G_2/T^2$,
и $W$ --- его группа Вейля. действующая на $M$ справа,
$\sigma \in W$ --- отражение системы корней типа $A_3$ или $G_2$.

\begin{PROP}{}\label{PROP:1.3} Для
флаговых пространств $M=SU_4/T^3$ и $G_2/T^2$
соответственно имеем $\Cal E(M/\sigma ) = \nu (M/\sigma ) = 12$ и $16$.
\end{PROP}

Пусть $M=G/H$ -- флаговое пространство (ф.п.), $\Omega $ --- его система $T$-корней,
$W=\opn{Norm}_G(H)/H$ -- его группа Вейля, действующая на $M$ справа,
$\sigma \in W $  -- элемент, индуцирующий отражение системы $\Omega $,
и $M/ \sigma $ -- факторпространство по группе $\langle \sigma  \rangle$.

\begin{PROP}{}\label{PROP:1.4} Если $\Omega $ является системой корней типа $G_2$, то
$$
\Cal E(M/\sigma )= \nu(M/\sigma ) = 16.
$$
\end{PROP}

Предложения~\ref{PROP:1.3} и~\ref{PROP:1.4} следуют из лемм,
которые будут доказаны
в пп.~\ref{sect:8.1}--\ref{sect:8.3}.

Существует пять флаговых пространств $M$ с системой $T$-корней типа $G_2:$
\begin{equation}{}\label{eq:15}
G_2/T^2, F_4/T^2A_2, E_6/T^2(A_2)^2, E_7/T^2A_5, E_8/T^2E_6,
\end{equation}
где подгруппа $A_5$ группы $E_7$ не содержится в подгруппах типа $A_6$.
%
%
Каждому из этих 
ф.п. соответственно
можно сопоставить единственную схему Дынкина:
$$
\let\0\circ
\let\1\bullet
\begin{smallmatrix} \1 \,\Lleftarrow \,\1 \end{smallmatrix}\,, \quad
\begin{smallmatrix} \0&\0\, \Leftarrow \,\1&\1 \end{smallmatrix}\,, \quad
\begin{smallmatrix} \0&\0&\1&\0&\0\\&&\1 \end{smallmatrix}, \quad
\begin{smallmatrix} \0&\0&\0&\0&\1&\1\\&&&\0 \end{smallmatrix}, \quad
\begin{smallmatrix} \1&\1&\0&\0&\0&\0&\0\\&&&&\0 \end{smallmatrix} 
.
$$
Следовательно, каждое из них имеет единственную, с точностью до изоморфизма,
комплексную форму $C$.
Как и вообще для ф.п. \cite{Ale-P},
формы $C$ однозначно соответствуют камерам системы $\Omega $.
Поскольку $\Omega $ является системой корней,
эти камеры совпадают с классическими камерами Вейля.
Следовательно, группа $W$ ин\-ду\-цирует в $\RR^2 = \RR \Omega $
группу линейных преобразований, сохраняющую систему кор\-ней  типа $G_2$
и транзитивную на множестве камер Вейля;
т.е. всю группу Вейля этой системы корней.
Это доказывает существование $\sigma $
для каждого из пространств \eqref{eq:15}.



Многогранник Ньютона $P=P(G_2)$ пространств \eqref{eq:15}
является одним из многогранников $P(\Omega )$,
где $\Omega $ --- система корней,
заданных общей формулой \cite[(0.1)]{2007}.
Он будет изучен в пп.~\ref{sect:8.4}--\ref{sect:8.5}.

\begin{PROP}{}
Для каждого флагового пространства $M$ вида \eqref{eq:15}
число Ньютона $\nu(M)$ выражается формулой
\begin{align}{}
\label{eq:nu=152} \nu (M) &= 152.
\end{align}
\end{PROP}

Доказательство приводится в п.~\ref{sect:8.4}.

В п.~\ref{sect:8.6} доказано неравенство $\Cal E(G_2/T^2)< \nu (G_2/T^2)$.

Дефект $ \delta _M = \nu (M)- \Cal E(M)$ можно найти
независимо от \eqref{eq:nu=152}.

\begin{PROP}{}\label{PROP:8.4}
Для каждого флагового пространства $M$ вида \eqref{eq:15} имеем:
\begin{align}{}
\label{eq:152-0}  \nu(M) &- \Cal E(M) = 0   \quad\text{при}\quad M\ne G_2/T^2,\\
\label{eq:152-18} \nu(M) & - \Cal E(M) = 18 \quad\text{при}\quad M=G_2/T^2.
\end{align}
Пространство $M=G_2/T^2$ обладает положительным дефектом $\nu(M) - \Cal E(M)$,
накопленным в $12$ сжатых однородных пространствам $M_{\gamma }$,
соответствующих квадратным граням $\gamma $ многогранника Ньютона $P$.
\end{PROP}

Хотя доказательство предложения~\ref{PROP:8.4} 
в принципе не сложно и не громоздко,
объем настоящей работы все же не позволяет его включить.



В следующих замечаниях $M$ --- флаговое пространство вида \eqref{eq:15}.

\begin{REM}{} Число г.к. положительно определенных инвариантных метрик Эйнштейна
в $M/\sigma $ равно $2$ при $M=G_2/T^2$ и $6$ в остальных четырех случаях
(см. добавление).  Обратные образы этих метрик заведомо не являются метриками
Кэлера--Эйнштейна по отношению к какой либо инвариантной комплексной
структуре $J$ на $M$.
\end{REM}

\begin{REM}{} Комплексные метрики Эйнштейна на $M/\sigma $ смотри
в добавлении~2 в конце раздела. Из их обратных образов действием
группы $W$ получается $2+3 \times 14 = 44$ г.к. метрик на $M$.
Присоединяя метрики Кэлера-Эйнштейна относительно различных $J$,
получаем $50$ г.к. комплексных инвариантных метрик Эйнштейна на $M$.

Пусть $k$ --- число остальных таких метрик Эйнштейна на $M$
с точностью до гомотетий и до действия группы $W$.
Приняв гипотезу об однократности всех решений,
нетрудно доказать, что $6k=\nu(M)- 50 - \delta _M $.
Тогда в силу \eqref{eq:nu=152} и \eqref{eq:152-18}
$$
k = \frac{1}{6}\,(152 - 50 - 18) =14
\qquad \text{при $M=G_2/T^2$}.
$$
\end{REM}

\begin{REM}{}\label{REM:8.3}
Для доказательства предложения~\ref{PROP:8.4}
надо прежде всего описать грани многогранника Ньютона $P$.
Можно показать, что $P$ имеет $15$ двумерных квадратных граней,
три из которых отвергаются
по признаку \cite[Test 7.2 (\S7.1)]{2007}.
Все грани в форме пирамиды можно отвергнуть
по признаку \cite[Test 7.1]{2007}.
Анализ
оставшихся граней сводится к $8$ случаям,
каждый из которых можно разобрать коротко и элементарно.
В результате все не квадратные грани отвергаются.

Для доказательства равенства \eqref{eq:152-18}
рассматривалось торическое многообразие $\widetilde V$
(см. \cite[замечание 3.5]{2007})
и некоторые из его однооорбитных стра\-тов,
а именно, (неособые) стра\-ты,
соответствующие $12$ квадратным граням.
Оценка $\delta _M\ge 18$ при $M=G_2/T^2$
следует из простых и достаточно общих соображений, а все необходимое для доказательства
равенства $\delta _M=18$ (после некоторой подготовки) обнаруживается на одном
развороте книги \cite[гл. 1, \S 9.7, с.128]{AVGU}.
\end{REM}



\subsection{Инварианты де Зибенталя пространств $M/\sigma $}
\label{sect:8.1}
(См. \cite[Proposition 2.4, Example 1.1]{2007}).
Найдем инварианты де Зибенталя
однородных пространств $M/ \sigma $,
где $M=SU_4/T^3$
или $M$ --- любое из пяти флаговых пространств \eqref{eq:15}. Отсюда
получим соответствующие многогранники Ньютона и числа Ньютона.

Фиксируем тетраэдр $\Delta _3$, вершину $\beta \in \Delta _3$
и симметрию $\tau \in {\mathrm O}(3): \Delta _3 \to \Delta _3 $
такую, что $\tau (\beta ) = \beta $ и $\tau ^2 = e$.
Очевидно, $\tau $ --- отражение.

Определим тройные отношения $\opn TA_3$ и $\opn TG_2$
на множестве $I$ неориентированных ребер тетраэдра условиями
''ребра образуют треугольник'' и ''ребра образуют треугольник или
три ребра пересекаются в точке $\beta $''.

Пусть $I=\{\gamma _1,\gamma _2,\gamma _3 , \beta _1,\beta _2,\beta _3 \} $,
где ребра $\beta _i$ выходят из вершины $\beta $. Тогда
{\let\b\beta\let\g\gamma \begin{align}{}\label{eq:13}
\opn TA_3 &= \opn{Sym}\{
                  (\b_1,\b_2,\g_3),
                  (\b_1,\g_2,\b_3),
                  (\g_1,\b_2,\b_3),
                  (\g_1,\g_2,\g_3)\},
\\
\label{eq:14}
\opn TG_2 &= \opn{Sym}\{
                  (\b_1,\b_2,\g_3),
                  (\b_1,\g_2,\b_3),
                  (\g_1,\b_2,\b_3),
                  (\b_1,\b_2,\b_3), (\g_1,\g_2,\g_3) \},
\end{align}}%
где $\opn{Sym}$ означает переход к симметричному множеству троек.

\begin{LEM}{}\label{LEM:1.3}
Пусть $M=SU_4/T^3$ или $G_2/T^2$. Тогда
(a)
неприводимые компоненты представлений изотропии
пространства $M$ можно занумеровать шестью ребрами
$\beta _i,\gamma _i \in I$
так, что $\sigma \mid_{\,I} = \tau $, а
инвариант де Зибенталя для $M$
записывается в виде \eqref{eq:13} или \eqref{eq:14} соответственно;

(b)
инварианты де Зибенталя пространств $M/\sigma $     
имеют
соответственно 
вид:
\begin{align}{}\label{eq:11}
\opn T & = \opn{Sym} \{ (a',b',b), (a,b',b'),             (a,b,b) \}
,
\\              \label{eq:12}
\opn T &= \opn{Sym} \{ (a',b',b), (a,b',b'), (a',b',b'), (a,b,b) \}
.
\end{align}
%
%
\end{LEM}

\begin{proof}{}
Обозначим пары противоположных корней в системах $A_3$ и $G_2$ через
$\beta (i) = \{\pm(\varepsilon _i - \varepsilon _4)\}$,
где $i=1,2,3$, и
$\gamma (i) = \{\pm(\varepsilon _j - \varepsilon _k)\}$,
где $\{i,j,k \} = \{1,2,3 \}$.
Здесь $\varepsilon _1 +\varepsilon _2+\varepsilon _3+ \varepsilon _4 =0$.
Предполагается, что $\varepsilon _4=0$ в случае системы $G_2$.
Без потери общности,
$\sigma $ переставляет друг с другом $\beta  (1)$, $\beta (2)$ и $\beta (3)$.
Положив $\beta _i = \beta (i)$, $\gamma _i = \gamma (i)$, $i=1,2,3$,
получаем (a).
Выведем (b) из (a).
Пусть $\tau (\beta _1) {=} \beta _1$, $\tau (\gamma _1) {=} \gamma _1$,
$\tau (\beta _2) {=} \beta _3$, $\tau (\gamma _2) {=} \gamma _3$.
Обозначим через $\mathfrak m \subset \mathfrak g $ ортогональное
дополнение подалгебры изотропии относительно формы Киллинга.
Его разложение на не\-при\-водимые модули
группы изотропии однородного пространства $M/\sigma $ имеет вид
%
%
$$
\mathfrak m = \mathfrak m_{a} + \mathfrak m_{b} + \mathfrak m_{a'} + \mathfrak m_{b'}      \,\,,
$$
где
$\mathfrak m_{a}  = \mathfrak m_{\gamma(1)}$,
$\mathfrak m_{b}  = \mathfrak m_{\gamma(2)} + \mathfrak m_{\gamma(3)}$,
$\mathfrak m_{a'} = \mathfrak m_{\beta(1)}$,
$\mathfrak m_{b'} = \mathfrak m_{\beta(2)}+\mathfrak m_{\beta(3)}$.
Подставляя в \eqref{eq:13} и \eqref{eq:14}
${a} $ вместо ${\gamma(1)}$,
${b} $ вместо ${\gamma(2)}$ и  ${\gamma(3)}$,
${a'}$ вместо ${\beta(1)}$,
${b'}$ вместо ${\beta(2)}$ и ${\beta(3)}$,
получаем инварианты де Зибенталя $\opn T/\sigma $ пространств $M/\sigma $
в виде \eqref{eq:11} и \eqref{eq:12}:
{\let\b\beta\let\g\gamma \begin{align*}{}
\opn TA_3/\sigma  &= \opn{Sym}\{
                  (a',b',b),
                  (a',b,b'),
                  (a,b',b'),
                  (a,b,b)\},
\\
\opn TG_2/\sigma  &= \opn{Sym}\{
                  (a',b',b),
                  (a',b,b'),
                  (a,b',b'),
                  (a',b',b'), (a,b,b) \}.
\end{align*}}%
%
%
Это доказывает (b) леммы. 
\end{proof}

Аналогично имеем :

\begin{LEM}{} Пусть $M$ --- флаговое пространство с системой $T$-корней типа
$G_2$ (т.е. любое из пяти флаговых пространств \eqref{eq:15}).
Тогда инварианты де Зибенталя для $M$ и $M/ \sigma $
соответственно имеют вид \eqref{eq:14} и \eqref{eq:12}.
\end{LEM}


\subsection{Трехмерные многогранники $P_{12}$ и $P_{16}$}
\label{sect:8.2}

Опишем трехмерные многогранники $P_{12}$ и $P_{16}$ в $\RR^4$ ---
многогранники Ньютона однородных пространств $M/\sigma $.
На этот раз для удобства обозначим стандартный базис пространства
$\RR^4$ через $\{a,a',b,b' \}$.
Рассмотрим следующие семь точек в $\RR^4:$
$$
\begin{array}{llll}
A=-a,\quad   &  B=b'-b-a',\quad   &  C=b-a'-b' \quad & \\
D=a'-b-b' , &E=a-2b, & F=a-2b' & G= a'-2b'
\end{array}
$$
\par

Обозначим через $P_{12}$ призму с треугольными основаниями
$(A,B,C)$ и $(D,E,F)$. Постоим над ее квадратной 
гранью $(D,A,C,F)$ пирамиду $P_{4}$ с вершиной $G \notin P_{12}$.

\par

Обозначим через $P_{16} = P_{12}\cup P_{4} $ выпуклую оболочку
всех семи точек $A, \dots ,G$ (Рис.~3). Многогранник $P_{16}$
комбинаторно устроен так же, как выпуклая оболочка семи вершин куба;
он имеет семь вершин и семь двумерных граней, а именно, три четырехугольные
и четыре треугольные грани:
$$
(A,B,E,D),(D,E,F,G),(B,C,F,E), 
(A,B,C),(A,C,G),(A,D,G),(C,F,G).
$$

Из определений следует:

\begin{LEM}{} Построенные многогранники $P_{12}$ и $P_{16}$
являются многогранниками Ньютона
однородных пространств класса $\nnnn$
(см. \cite[\S 2.4]{2007})
с инвариантами де Зибенталя \eqref{eq:11} и \eqref{eq:12} соответственно.
\end{LEM}

\begin{proof}{} Лемма следует из включений
$
S \subset P_{12} \subset P_{16}, 
$
где $S$ --- стандартный симплекс в $\RR^4$, т.е. симплекс с вершинами
$-a,-a',-b,-b'$.
\end{proof}

\begin{COR}{}
$P_{M/\sigma } = P_{12}$ для $M=SU_4/T^3$; \  \
$P_{M/\sigma } =P_{16}$ для $M=G_2/T^2$
и остальных флаговых пространств \eqref{eq:15}.
\end{COR}

\begin{LEM}{}\label{LEM:1.1} Приведенные объемы $\nu(P)$
многогранников $P_{12}$ и $P_{16}$ равны $12$ и $16$.
\end{LEM}

\begin{proof}{} Натянем на многогранники трехмерное аффинное пространство
$D = \{x \in \RR^4: \sum x_i = -1 \} $.
Существует куб $Q$ в $D$ со стороной евклидовой длины $2$ такой, что
$Q \supset P_{12}$.
Имеем $\nu(Q) = 3\cdot 2^3 =2\nu(P_{12})$, откуда $\nu(P_{12}) = 12$.
Приведенный объем трехмерной пирамиды в $D$ выражается необычной формулой
$$
\nu = h\, \sigma ,
$$
где $h$ --- высота пирамиды и $\sigma $ --- евклидова площадь основания.
Для пирамиды $P_{4}$ имеем $\sigma = 2 \times 2 $ и $h=1$, поскольку
$2h$ --- диагональ квадрата $(D,G,F,\frac{E+F}{2})$ со стороной $\sqrt{2}$.
Значит, $\nu (P_{16})= \nu(P_{12}) + \nu(P_{4}) = 12+4$.
\end{proof}

\begin{figure}{}
\phantom{.} \par
\parbox{13cm}{\includegraphics{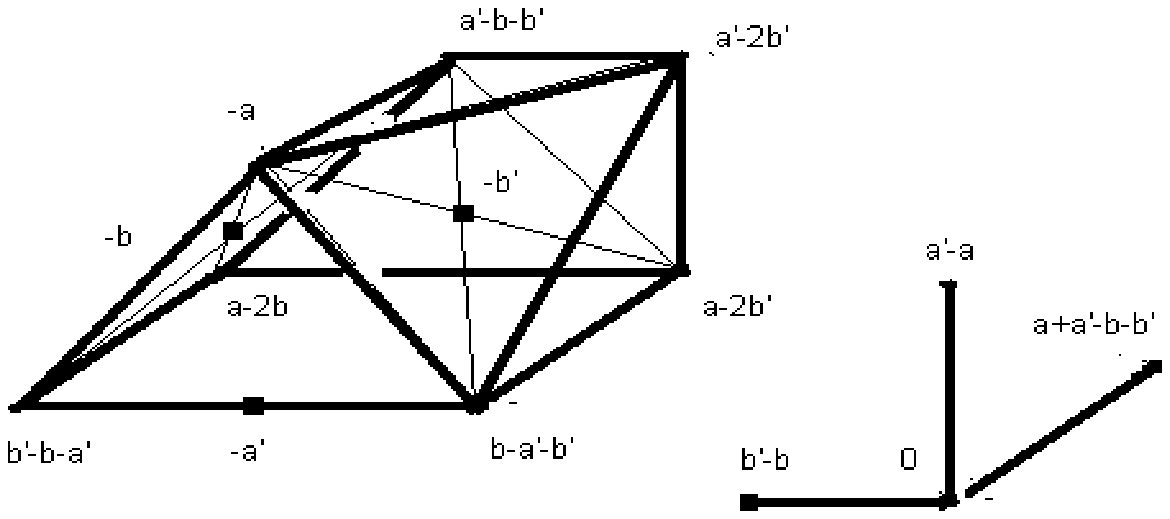}}
\par
Рис.3. Многогранник $P_{16}$ и естественный
\\
ортогональный репер в точке $-a'\in P_{16}$.
\end{figure}



\subsection{Однородные пространства с трехмерными многогранниками $P_{12}$ и $P_{16}$}
\label{sect:8.3}

\begin{LEM}{}\label{LEM:1.2} Пусть $M=G/H$ --- компактное
однородное пространство с однократным  спектром представления изотропии,
с многогранником Ньютона $P=P_{12}$ или $P_{16}$.
Тогда группа $G$ полупроста, $M$ содержится в классе $\nnnn$,
определенном в
\cite[\S 2.4]{2007}.
Равенство
$$
\Cal E(M) = \nu (P)  = 12, \,16
$$
эквивалентно
$[\mathfrak m_{a'},\mathfrak m_{a'}]\ne0$
и следует из
$[\mathfrak m_{a'},\mathfrak h]\ne 0$.
\end{LEM}

\begin{proof}{} Из включений $S \subset P_{12} \subset P_{16}$ следует,
что группа Ли $G$ полупроста
(в силу \cite[\S1.6, замечание 1.1, 2)]{2007})
и, более того, $M$ содержится в классе $\nnnn$
(в силу \cite[предложение 2.2]{2007}).
Пусть $t=(t_{a},t_{b},t_{a'},t_{b'})$ ---
естественные координаты инвариантной комплексной метрики в $M$
такие, как в \cite[пример 1.1]{2007},
а $s(t)$ --- скалярная кривизна.
Носитель $\opn{supp}(s)$ лорановского многочлена $s(t)$
состоит из вершин многогранника
$P$ и точек $-a',-b,-b' \in P$ 
(\cite[предложение 2.4]{2007}).

Существует единственное ребро $\gamma $ с центром в точке носителя.
Это ребро $\gamma = (BC)$ с центром $-a' \in \opn{supp(s)} $
\textit{отвергается}
(в смысле \cite[\S7.1]{2007})
при $[\mathfrak m_{a'},\mathfrak m_{a'}]\ne0$
по признаку \cite[Test 7.2]{2007}.

В соответствии 
с \cite[\S 7.1]{2007},
проверим, что
отвергаются все остальные $p$-грани $\gamma \subset P$, ${ 0 < p < \dim(P)}$,
и тем самым докажем лемму.
Как видно из чертежа,
все треугольные грани $\gamma $ и все ребра $\gamma \ne (BC) $
отвергаются
по признаку \cite[Test 7.1]{2007}.
Осталось рассмотреть три четырехугольные грани $\gamma \subset  P$,
где $P=P_{12}$ или $P_{16}$.

Случай 1). Пусть $\gamma = (A,B,E,D)$ -- квадратная грань, $P=P_{12}$ или $P_{16}$.
Тогда
$$
\phi (t) := t_{b} \,s_{\gamma }(t)
=   A\,t^{b-a} - E\,t^{a-b}   -  B\,t^{b'-a'} - D\,t^{a'-b'} + K, \qquad B=D.
$$
(на этот раз мы обозначаем через, $A$, $B$, и т.д. не вершину многогранника
$P$, а соответствующий коэффициент многочлена $\pm s(t)$).
Функция $\phi $ существенно зависит от $2=\dim(\gamma )$ переменных.
Переходя от четырех переменных $t_i$ к двум новым переменным $x,y$, получаем
$$
\phi (x,y) = Ax^{-1}-Ex - B(y^{-1}+y) + K.
$$

Случай 2). Пусть $\gamma = (B,C,E,F)$ -- прямоугольная грань, $P=P_{12}$ или $P_{16}$.
Тогда для двух подходящих переменных $x,y$
$$
\phi (t) := - t_{a'} \,s_{\gamma }(t) =
 By^{-1} + Cy + Exy^{-1} + Fxy - L, \qquad B=C.
$$

Случай 3). Пусть $\gamma = (D,E,F,G)$ -- трапеция, $P=P_{16}$.
Тогда для двух подходящих переменных $x,y$
$$
\phi (t) :=  t^{-(a-b-b')} \,s_{\gamma }(t) =
- Ey^{-1} - Fy - Dx - Gxy.
$$

Покажем, что в каждом из этих трех случаев комплексная 
кривая $Y_{\gamma }$ в ${(\CC \setminus 0)^2}$,
заданная уравнением $\phi (x,y) = 0$, $x,y\ne0$, является неособой.
Пусть
$$
z = \partial \phi (x,y)/\partial x, \qquad
\psi =  \phi -xz.
$$
В случае 1) имеем $z= -Ax^{-2}-E$, т.е. $x^{-2} = - (z+E)/A$.
Из $\partial \phi (x,y)/ \partial y =0  $ следует $y^2=1$.
Тогда $\psi (x,y) = 2Ax^{-1} \pm 2B +K$,
\begin{equation*}{}
l_1 :=
\psi (x,y)\,\psi (-x,y)
= (K\pm 2B)^2 - 4A^2x^{-2} = (K\pm 2B)^2 + 4AE + 4Az.
\end{equation*}
В случае 2) имеем $z=Ey^{-1}+Fy$, \enskip  $\psi (x,y) = By^{-1}+By -L$.
Значит, $Ey^{-2}= - F+y^{-1}z$, $Fy^{2}= - E+yz$,
\begin{equation*}{}
\begin{aligned}{}
l_2 :=
\psi (x,y)\,\psi (x,-y) &= L^2 - B^2(y^{-2} + y^{2} +2)
        = L^2 + B^2(F/E+E/F-2)
        \\ &
        - (B^2/EF)(Fy^{-1}+Ey)\,z
\end{aligned}
\end{equation*}
В случае 3) имеем $ z=-D-Gy$, т.е. $ y = -(D+z)/G$,
$y^{-1} =  -(y^{-1}z+G)/D$,
\begin{equation*}{}
l_3 :=
\psi (x,y) = -Ey^{-1}-Fy = (GE/D+DF/G) + ((E/D)y^{-1}+F/G)\,z .
\quad
\end{equation*}
Здесь $A,K,L$ --- коэффициенты многочлена $s(t)$
соответственно при $t_{a}^{-1}=t^{-a}$,
$t_{b}^{-1}=t^{-b}$ и $t_{a'}^{-1}=t^{-a'}$.
%
%
Из компактности $M$ следует,
что они неотрицательны 
(\cite[лемма 5.3, 2)]{2007})
(срав. также с формулой Вана--Циллера в \cite[пример 1.1]{2007}).
Тогда ${A,K,L>0}$ в силу
\cite[предложение 2.4]{2007} о носителе многочлена $s(t)$.
Имеем:
$$
A,B,C,D,E,F,G,K,L > 0, 
$$
следовательно, выполняются дискриминантные неравенства
$$
(K\pm 2B)^2 + 4AE>0,\qquad
EFL^2 + B^2(E-F)^2>0,\qquad
EG^2 + FD^2>0,
$$
в силу которых $l_1,l_2,l_3>0$ при $z=0$.
Поэтому кривые $Y_{\gamma }$ не имеют особых точек, и
грани $\gamma = (A,B,E,D)$, $(B,C,E,F)$, $(D,E,F,G)$ отвергаются.

Случай 3a). Грань $(D,A,C,F)$ призмы $P_{12}$ (содержащая внутреннюю точку $-b'$ многогранника $P_{16}$)
рассматривается точно так же, как грань $(A,B,E,D)$ в  $P_{12}$,
и так же отвергается. Лемма доказана.
\end{proof}



\begin{REM}{}
Для того, чтобы однородное пространство
$M = G/H$ компактной
полупростой группы Ли $G$ с инвариантом 
\eqref{eq:11} или \eqref{eq:12} соответственно
имело многогранник Ньютона $P=P_{12}$ или $P_{16}$, достаточно,
чтобы выполнялось $[\mathfrak m_{a},\mathfrak h]\ne 0$.
\end{REM}

Приведем другие примеры однородных
пространств с инвариантом \eqref{eq:11}.

\begin{EXAM}
\label{EXAM:8.1}
Инвариантом де Зибенталя следующих 
однородных пространств является тройное отношение
\eqref{eq:11}:
\begin{enumerate}
\item
односвязные пространства
\footnote{
В \cite[\S 7]{2007}, 
подстрочное примечание 11),
одно пространство из следующего списка было пропущено
($E_7/T^1{\cdot}D_4{\cdot}A_1{\cdot}A_1$)
и приклеились два пространства из другого списка
(а именно, особые флаговые пространства с системой $T$-корней
$\Omega = \{\pm1,\pm2,\pm3 \} \subset \ZZ^1$, их многогранник Ньютона --- трапеция).
}
\begin{center}{}
$E_7/T^1{\cdot}A_1{\cdot}A_1{\cdot}D_4$,\quad
$E_8/T^1{\cdot}A_1{\cdot}D_6$,\quad  
$F_4/T^1{\cdot}A_1{\cdot}B_2$;
\end{center}
%
%
%
%
\item
три серии
(где $l{=}m'{+}n{+}m$):
\begin{center}{}
$
\Sp_{\,l}/\Sp_{\,m'}{\times} \oU_n {\times} \Sp_{m}\,( {1{\le} m{\le} m'} ),
$
\,
$
\SO_{\,2l}/\SO_{2m'}{\times} \oU_n {\times} \SO_{2m}\,( {2{\le} m{\le} m'} ),
$
\\
$
\SO_{\,2l+1}/\SO_{2m'}{\times} \oU_n {\times} \SO_{2m+1}\,( { 0{\le} m,\,2{\le} m' } ),
$
\end{center}
\item
серия
пространств $G/H$ с $rank(G)>rank(H):$
\begin{center}{}
$G/H=\SU_n/\oS(\Sp_{\,n_1}{\cdot}T^1{ \times}\oU_{n_2}{ \times}\oU_{n_3})$
($n=2n_1+n_2+n_3$). 
\end{center}
\end{enumerate}
(В пп.1) и 2) перечислены все $p$-симметрические пространства
внутреннего типа с системой $S$-корней $\ov{BC_2}$,
см. \cite[0.3.9]{2008}).
Пространства 1),2),3)
имеют многогранник Ньютона $P_{12}$. 
\end{EXAM}



Из лемм вытекает :

\begin{COR}{} $\Cal E(M/\sigma )= \nu (M/\sigma ) \in\{ 12,\, 16 \}$ для $M=SU_4/T^3$,
$G_2/T^2$ и других флаговых пространств \eqref{eq:15};
\\
$
\Cal E(M)= \nu (M) = 12
$
для однородных пространств $M$,
перечисленных в примере~\ref{EXAM:8.1}.
\end{COR}

\begin{REM}{} Призма $P_{12}=P(\ov{BC_2})$ отличается от
рассмотренной в
\cite[\S 7.2]{2007}
призмы $\mathrm{IB}=P(B_2)$.
Их можно перевести друг в друга действием группы $GL(4,\ZZ)$
(отсюда равенство объемов).
Но это не распространяется на носители $\opn{supp}(s)$.
\end{REM}



\subsection{Число Ньютона $\nu (G_2/T^2)$}
\label{sect:8.4}

Приведем доказательство формулы \eqref{eq:nu=152},
опуская только легко восстановимые подробности.
Пусть $P=P(G_2)$ --- пятимерный многогранник Ньютона
флаговых пространств $M$ с системой $T$-корней типа $G_2$.
По определению числа Ньютона имеем
$$
\nu  = \nu(M) = \frac{5!}{\sqrt{6}}\,\opn{vol}_E (P(G_2)). 
$$
Здесь 
$\opn{vol}_{E}(Q)$ --- стандартный евклидов объем любого многогранника
$Q$ в гиперплоскости $D=\{x \in \RR^6 : x_1+ \ldots +x_{6} = -1 \}$.
Сравнение формул \eqref{eq:13} и \eqref{eq:14} 
показывает, что 
$P=P(G_2)$ является выпуклой оболочкой рассмотренного ранее в \cite{2007}
многогранника Ньютона $\Pi = P(A_3)$ и трех точек. Обозначим эти
точки через $\alpha _i$, $i=1,2,3$. Отождествим
$I=\{\gamma _1,\gamma _2,\gamma _3 , \beta _1,\beta _2,\beta _3 \}$
со стандартным базисом
\footnote{
Здесь
от обозначений \cite{2007}
удобно перейти к новым обозначениям для базиса в $\RR^6$.
В формулах для $\Pi $
(см. \cite[пример 1.4]{2007})
можно положить
$\mathbf 1_{ij}=\mathbf 1_{ji}= \gamma _k$,  $\mathbf 1_{i4}= \beta _i$
при $\{i,j,k \}=\{1,2,3 \}$.
}
в $\RR^6$.
Тогда
$$
\alpha _i = \beta _i - \beta _j - \beta _k ,\quad \{i,j,k \}=\{1,2,3\},
$$
Обозначим через $\Cal C(Q )=\{q \in \RR^6 : \langle p,q \rangle\le0\,\,\forall\,p \in Q\}$
двойственный конус для $Q \subset  \RR^6$.
Пусть $\Gamma $ --- фасета многогранника $\Pi $
и $f \in \Cal C(\Pi ) \setminus 0$ --- вектор нормали к $\Gamma $
(параллельный ребру конуса $\Cal C(\Pi )$).
При $\langle \alpha _i,f  \rangle > 0$ назовем фасету $\Gamma $  {\bf освещенной}
из точки $\alpha _i$.
Из включений
$
\frac12\,(\alpha _i + \alpha _j) = - \beta  _k \in \Pi 
$, где $\{i,j,k \}=\{1,2,3 \}$,
легко следует, что
$P$ совпадает с объединением многогранника $\Pi $
и всех пирамид с вершиной $\alpha _i$,  построенных над фасетами, освещенными
из точки $\alpha _i$, $i=1,2,3$, причем эти многогранники образуют разбиение
многогранника $P$, т.е. имеют попарно непересекающиеся внутренности.


Пятимерный многогранник $\Pi = P(A_3)$
рассматривался в \cite[примеры 1.4, 3.2, 7.7]{2007}. Он
имеет семь фасет, т.е.
двойственный конус 
$\Cal C(\Pi )$ имеет семь ребер. Эти ребра параллельны следующим
характеристикам, т.е. $0,1$-векторам:
$$
\chi^0 := \chi_{\gamma _1,\gamma _2,\gamma _3},\qquad
\chi ^{i} = \chi _{\gamma _i,\beta _j,\beta _k},\qquad
\phi ^{i} = \chi _{\gamma _i,\beta _i} ,\qquad i=1,2,3.
$$
Здесь $\chi _S$ --- характеристика с нуль-множеством $S \subset I$,
$I=\{\gamma _1,\gamma _2,\gamma _3 , \beta _1,\beta _2,\beta _3 \}$,
т.е. вектор а координатами $\chi_S(i)=0$ при $i \in S$, $1$ при $i\in I \setminus  S$.
Для других $G/H$ характеристики $\chi_S$ определяются аналогично.
Вос\-пользуемся тем, что в общем случае $0,1$-векторы $\chi _S \in  \Cal C(P_{G/H})$,
и только они, являются характеристиками $H$-инвариантных подалгебр
$\mathfrak k = \mathfrak h + \mathfrak m_S \subset \mathfrak g$.
(ср. \cite[следствие 2.0]{2008} или \cite[\S 6.1.5]{2007}).
Например, при $G/H{=}SU_4/T^3$ имеем $\chi^0$, $\chi^i$, $\phi ^i \in \Cal C(\Pi )$,
поскольку $P_{G/H}{=}\Pi$ -- многогранник Ньютона пространства $SU_4/T^3$,
а $\chi^0$, $\chi^i$, $\phi ^i$ -- характеристики максимальных подалгебр максимального ранга
алгебры Ли $\mathfrak {su}_4$.
При $G/H{=}G_2/T^2$ имеем $P_{G/H}{=}P$.
Легко видеть, что $\chi^0$, $\phi ^i$ -- характеристики подалгебр $\mathfrak k$
особой алгебры Ли $\mathfrak g_2$.
(Структура всех подалгебр алгебры $\mathfrak g_2 $, собственных и
собственным образом содержащих $\mathfrak h= Lie(T^2)$, является де\-ре\-вом,
изображенным на Рис.4, где
черными кружками обозначены подалгебры $\mathfrak k$ с характеристиками
$\chi^0$ и $\phi ^{i}$, $i=1,2,3$,
белыми кружками --- другие подалгебры $\mathfrak k$,
отрезками --- включения $\subset $ и $\supset $ между подалгебрами).
\begin{figure}{}
\begin{center}
\parbox{6cm}{\includegraphics{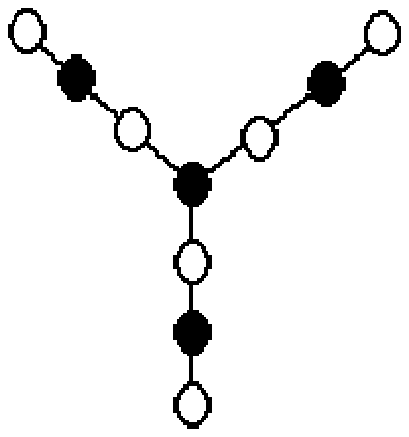}}
\parbox{4cm}{Рис.4. Граф Вана--Циллера пространства $G_2/T^2$.\\
(Вершины графа отвечают подалгебрам $\mathfrak k$, $ \mathfrak t_2 \subsetneqq  \mathfrak k \subsetneqq \mathfrak {g}_2$)}
\end{center}
\end{figure}
Следовательно, $\chi^0$ и $\phi ^{i}\in \Cal C(P) \subset \Cal C(\Pi )  $.
Значит, фасеты многогранника $\Pi$
с векторами нормалей $\chi^0$ и $\phi ^{i}$, $i=1,2,3$, лежат на фасетах
многогранника $P$, и потому не освещены.
Далее, имеем:
$$
\langle \alpha _i,\chi ^{\,j} \rangle = 2 \delta _{i}^j-1, \quad i,j=1,2,3,
$$
где $\delta _{i}^j$ --- символ Кронекера.
Следовательно, $i$-я фасета
$
\Gamma _i = \{p \in \Pi  : \langle p,\chi ^{\,i} \rangle = 0 \}
$
многогранника $\Pi $
освещена из $i$-й точки $\alpha _i$, $i=1,2,3$. Обозначим
через $P_i$ пирамиду с апексом $\alpha _i$, построенную на фасете $\Gamma _i$.
Итак, мы разбили многогранник $P=P(G_2)$ на $4$ многогранника
с попарно не пересекающимися внутренностями:
$$
P(G_2) = P_1 \cup P_2 \cup P_3 \cup \Pi ,\qquad \Pi := P(A_3).
$$
Нетрудно проверить, что высота $h$ каждой пирамиды $P_i$ равна ${\sqrt{6}}/{3}$.

Подсчитаем объем $P(G_2)$.
Многогранник $\Pi $ является евклидовым прямым произведением
тетраэдра с ребром $2$ и правильного треугольника $\triangle_2$
со стороной $2$
(см. \cite[пример 1.4]{2007})
$\Pi = \triangle_3 \times \triangle_2$,
откуда, очевидно, $\frac{5!}{\sqrt{6}}\,\opn{vol}_E(\Pi ) = 80$
(\cite[пример 1.5]{2007}).
При этом $\Gamma _i = \triangle _2 \times \triangle _2$.
Поэтому евклидова площадь основания $\Gamma _i $ пирамиды $P_i$ равна $3$.
Отсюда находим объемы пирамид:
$$
\opn{vol}_{E}(P_i) = \frac{3h}{\dim(P_i)} = \frac{\sqrt{6}}{5},
\qquad i=1,2,3.
$$
Значит
$$
\nu  
= 3\cdot 4! + \frac{5!}{\sqrt{6}}\,\opn{vol}_E(\Pi ) = 72 + 80 = 152.
$$
Формула \eqref{eq:nu=152} доказана.



\subsection{Фасеты многогранника Ньютона $P=P(G_2)$}
\label{sect:8.5}
Продолжим изучение вы\-пу\-к\-лого пятимерного многогранника $P=P(G_2)$.
Найдем его фасеты $\Gamma $, 
проходящие через
точку $\alpha _1 = \beta _1 - \beta _2 - \beta _3 = (0,0,0,1,-1,-1)$.
Очевидно, $\Gamma$ содержит боковую $4$-грань пирамиды $P_1 = \opn{Conv}(\,\alpha _1,\,\Gamma _1\,)$.
%
%
Используя равенства
$$
\langle \alpha _i, \phi ^{\,j} \rangle = - 2 \delta\, _i^j,\qquad
\langle \alpha _i, \chi ^{\,j} \rangle\, = 2 \delta\, _i^j - 1, \qquad i,j=1,2,3,
$$
находим всего $6$ кандидатур 
для векторов нормалей $f \in \Cal C(P)$ к фасетам $\Gamma \subset P: $
\begin{equation*}{}
f \in \{ \phi ^2, \phi^3, \chi ^1 + \chi ^2, \chi ^1 + \chi ^3 , \chi ^0 + \chi ^1 , \phi ^1 + 2 \chi ^1  \}.
\end{equation*}
Все они будут векторами нормалей к фасетам
многогранника $P $,
поскольку пирамида $P_1$ заведомо не является $5$-симплексом
(напомним, что $\Gamma _1 = \triangle_2 \times \triangle_2 $).
%
%
%

\begin{PROP}{}\label{PROP:f} Конус $\Cal C(P(G_2))$ имеет $13$ ребер.
Пять из них порождены векторами\\
$ f_1 = (0,0,0,1,1,1) = \chi ^0 = \chi _{\gamma _1,\gamma _2,\gamma _3}, $\\
$ f_2 = (0,1,1,0,1,1) = \phi ^1 = \chi _{\gamma _1,\beta 1}, $\\
$ f_3 = (2,1,1,0,1,1) = \chi ^2 + \chi ^3, $ \\
$ f_4 = (0,1,1,2,1,1) = \chi ^0 + \chi ^1 = 2 - f_3, $ \\
$ f_5 = (0,3,3,2,1,1) = \phi ^1 + 2\chi ^1, $ \\
а остальные ребра получаются из ребер $\RR_+ f_i$, $i=2,3,4,5,$
действием шести\-чле\-н\-ной группы линейных пре\-образований
$\gamma  _i \mapsto \gamma _{\sigma (i)},\,\,
\beta _i \mapsto \beta _{\sigma (i)},\,\,   i=1,2,3$, 
где $\sigma  \in S_3$.
\end{PROP}

\begin{proof}{} Пусть $\Gamma $ --- фасета многогранника $P$,
нормальная к ребру $\RR_+ f$ конуса $\Cal C(P)$.
Осталось рассмотреть случай, когда $\alpha _1,\alpha _2,\alpha _3 \notin \Gamma$.
В этом случае $ \dim \Gamma \cap \Pi  = 4$ и фасета $\Gamma \cap \Pi$
многогранника $\Pi $.
не является освещенной ни из одной  точки $\alpha _i$, $i=1,2,3$.
Тогда, с точностью до множителя, $f = \chi ^0$.
%
%
\end{proof}

\begin{COR}{} Пятимерный многогранник $P=P(G_2)$ в $\RR^6$
имеет $13$ фасет, т.е. четы\-рех\-мерных граней. Равенство
$f_4 = 2 - f_3$, т.е. $f_4(i) = 2-f_3(i)$, $i=1, \dots ,6$,
показывает, что фасеты $\Gamma _{f_i}$, $i=3,4$, параллельны друг другу.
\end{COR}



\subsection{Неравенство для $\Cal E(G_2/T^2)$}
\label{sect:8.6}

Докажем, что $\Cal E(G_2/T^2)<\nu(G_2/T^2)$.
Вначале за\-ме\-тим, что
ана\-ло\-гичное неравенство для $SU_4/T^3$ можно доказать,
используя любую из $12$ двумерных квадратных граней $\square$ многогранника
Ньютона $\Pi = \triangle_3 \times \triangle_2$,
описанных
в \cite[\S 7.3]{2007},
и приравнивая нулю определитель второго порядка из~7.3.4.

Перенесем этот план на случай $G_2/T^2$.
Перейдем от $\Pi = P(A_3)$ к многограннику Ньютона $P=P(G_2)$
пространства $G_2/T^2$ и воспользуемся включением $\Pi \subset P$.
Утверждается, что все $12$ квадратов $\square \subset \Pi  $ являются
гранями многогранника $P$.
Для наших целей достаточно построить любой из них,
но мы построим два.
\footnote{
Попутно отметим, что
любую двумерную нетреугольную грань многогранника $P$
мож\-но перевести в одну и только одну из указанных ниже граней \emph{1--3}
посредством группы $S_3 = W_{G_2}/(\pm1)$.
}
В силу предложения~\ref{PROP:f} двойственный
конус 
$\Cal C(P)$ содержит следующие три вектора:
$$
\begin{array}{rrr}
{(1,1,2,4,2,3)}
\\ = (0,0,0,1,1,1) \\ +  (1,0,1,1,0,1)\\ + (0,1,1,2,1,1)\\  \phantom{()}
\end{array}
\qquad
\begin{array}{rrr}
{(1,5,6,3,3,2)}
\\ = (0,1,1,0,1,1)\\ + (1,1,2,1,1,0)\\ +  (0,3,3,2,1,1)\\ \phantom{()}
\end{array}
\qquad
\begin{array}{rrr}
{(0,5,5,5,4,4)}
\\ = (0,0,0,1,1,1)\\+ (0,1,1,0,1,1)\\ + (0,1,1,2,1,1)\\ + (0,3,3,2,1,1)
\end{array}
$$
Значит, ортогональное дополнение каждого из них пересекает
многогранник 
$P$ по некоторой грани.
Построенные грани лежат в $\Pi $, 
поскольку
не со\-дер\-жат вершин $\alpha _i \in P$, $i=1,2,3$,
вида $\alpha _i= \beta _i - \beta _j - \beta _k$,
на\-пример, вершины $\alpha _1 = {(0,0,0,1,-1,-1)}$.
Легко проверяется, что
\begin{enumerate}
\item[1.] $\Gamma _{(1,1,2,4,2,3)}$ --- квадрат вида $\square$,
\item[2.] $\Gamma _{(1,5,6,3,3,2)}$ --- квадрат вида $\square$,
\item[3.] $\Gamma _{(0,5,5,5,4,4)}=\Gamma _{(0,1,1,1,1,1)}$ --- квадрат
с центром 
$- \gamma _1 = (-1,0,0,0,0,0)$,\\
удовлетворяющий
признаку \cite[Test 7.2]{2007} (см. \cite[\S7.1]{2007}).
\end{enumerate}
Условие совместности системы уравнений 
\cite[(1.2)]{2007},
связанной с гранью $\square$, имеет вид
$$
D=\det\left\|\,\begin{matrix} a&b\\ c&d\end{matrix}\,\right\| =0 ,
$$
где $a,b,c,d$ --- коэффициенты многочлена Лорана $s(t)$
при степенях $t$, соответствующих вершинам грани $\square$.
Из формул для $s(t)$, приведенных в следующем дополнении, вы\-те\-кает, что
в случае $G_2/T^2$ коэффициенты при всех вершинах
многогранника $\Pi $ (но не $P$) совпадают, откуда $a=b=c=d$;
\enskip
$D=0$ в случае $G_2/T^2$ и $D\ne 0$ для
остальных флаговых пространств \eqref{eq:15}.
Это доказывает неравенство $\Cal E(G_2/T^2)<\nu(G_2/T^2)$.





\section*{Дополнение 1. О пяти особых флаговых пространствах $M$}


Рассмотрим группы изотропии $H$ и запишем многочлены Лорана $s(t)$ для пяти
флаговых пространств $M=G/H$ вида \eqref{eq:15}.
Тогда $H$ содержится в связной максимальной подгруппе максимального
ранга $A_2\cdot H'$, где $H'=(H,H)$ --- ее коммутант,
$(A_2,H')=\{e\,\}$.
Поэтому через точку $eH$ в $G/H$ проходит
$H$-инвариантное подпространство $ \CP^2 = A_2/T^2$.
В нем лежат три $H$-инвариантных подпространства $\CP^1$,
$eH \in \CP^1 \subset \CP^2$.
Их касательные пространства
$\mathfrak m_{\gamma (1)}$, $\mathfrak m_{\gamma (2)}$, $\mathfrak m_{\gamma (3)} \simeq\RR^2 $
являются тремя из шести неприводимых подмодулей $H$-модуля
$$
\mathfrak m = \mathfrak m_{\gamma (1)} + \mathfrak m_{\gamma (2)} + \mathfrak m_{\gamma (3)}
+ \mathfrak m_{\beta  (1)} + \mathfrak m_{\beta (2)} + \mathfrak m_{\beta (3)}
$$
(см. обозначения в доказательстве леммы~\ref{LEM:1.3}, ср. \cite[пример 1.1]{2007}).
Подмодули $\mathfrak m_{\gamma (i)}$ отвечают длинным корням системы $\Omega $
(которая является системой корней типа $G_2$).
Учитывая симметрии, находим размерности $N_{\alpha } = \dim (\mathfrak m_{\alpha }):$
$$
N_{\gamma (1)} = N_{\gamma (2)} = N_{\gamma (3)} = 2,
\qquad
N_{\beta  (1)} = N_{\beta  (2)} = N_{\beta  (3)} = N ,
$$
где $N=(\dim M - 6)/3$.

Для коммутирующих символов $x,y,z,x',y',z'$ введем обозначения
$$
((x,y,z)) = \frac{x^2+y^2+z^2}{xyz},
\quad
((x',y',z')(x,y,z)) = ((x',y,z)) + ((x,y',z)) + ((x,y,z')).
$$
Выражение для скалярной кривизны $s(t)$ из \cite[Example 1.1]{2007}
принимает вид 
\begin{equation}{}\label{eq:16}
\aligned
2s(t) &=
\sum_{i=1}^3 [2 t_{\gamma (i)}^{-1}+ N t_{\beta (i)}^{-1} ]
- \kappa_1 ((t_{\gamma (1)},t_{\gamma (2)},t_{\gamma (3)}))
\\    &
- \kappa_2 ((t_{\gamma (1)},t_{\gamma (2)},t_{\gamma (3)})(t_{\beta (1)},t_{\beta (2)},t_{\beta (3)}))
- \kappa_3 ((t_{\beta (1)},t_{\beta (2)},t_{\beta (3)})).
\endaligned
\end{equation}
Неравенства (Вана--Циллера) для коэффициентов 
из \cite[Example 1.1]{2007}
принимают вид
\begin{equation}{}\label{eq:17}
1- \kappa _1 - \kappa _2 > 0,     \qquad
N- 4 \kappa _2 - 2 \kappa _3 > 0, \qquad
\kappa _1,\kappa _2, \kappa _3 >0.
\end{equation}
Обозначим через $j$ индекс Дынкина присоединенного представления группы $G$.
Утверждается, что
\begin{equation}{}\label{eq:18}
N= j- 6,\qquad
\kappa _1 = 2/j,\qquad
\kappa _2 = 1-6/j,\qquad
\kappa _3 = j/6-1,
\end{equation}
при этом $j = 8,18,24,36,60 $ соответственно для $G=G_2,F_4,E_6,E_7,E_8$.




\section*{Дополнение 2. Вычисление комплексных метрик Эйнштейна в $M/\sigma $}

Пусть $M$ --- любое из флаговых пространств \eqref{eq:15},
$W$ --- его группа Вейля.


Нетрудно убедиться, что
на факторпространстве $M/W$ существуют две $G$-инвариантные
комплексные метрики Эйнштейна с точностью до гомотетии.
Их обратные образы на $M$ имеют координаты
$t_{\gamma (1)}=t_{\gamma (2)}=t_{\gamma (3)}$ и
\begin{equation}{}\label{eq:19}
t_{\beta (i)}/t_{\gamma (i)} =
\frac {1}{60}\,\left (5\,j\pm\sqrt {5
\left (j-12\right)\left(5\,j-12\right)}\right ), 
\quad i=1,2,3.
\end{equation}
%
%
%
При $M\ne G_2/T^2$ ($j\ge18$) формула \eqref{eq:19} задает две
не гомотетичные друг другу положительно определенные метрики Эйнштейна.

Перейдем к $G$-инвариантным метрикам в $M/\sigma $.
Пусть $t_{a}=t_{\gamma (1)}$, $t_{b}=t_{\gamma (2)}=t_{\gamma (3)}$,
$t_{a'}=t_{\beta (1)}$, $t_{b'}=t_{\beta (2)}=t_{\beta (3)}$
--- естественные координаты в пространстве этих метрик.

Следующая лемма показывает, что коэффициенты каждого не $W$-инвариантного
решения 
\cite[Eq. (1.1)]{2007}
(с точностью до общего множителя) принадлежат 
полю алгеб\-ра\-ических чисел степени $14 = \dim (G_2)$,
а поля для разных решений изоморфны.

\begin{LEM}{}\label{LEM:D2} На факторпространстве $M/\sigma $ существует $14$
попарно негомоте\-тич\-ных комплексных метрик Эйнштейна вида
$$
{{t_{{b'}}}/{t_{{b}}}}  = x, \quad
{{t_{{a'}}}/{t_{{a}}}} = F_1(x), \quad
{{t_{{b'}}}/{t_{{a}}}} = F_2(x), \quad  F_i(x) \in \QQ(x),
$$
где $x\in \{x_1, \dots ,x_{14} \} $
--- корень неприводимого многочлена $14$ степени $f_{14}(x) \in \QQ[x]$.
\end{LEM}

Лемма проверяется прямым вычислением с помощью программы MAPLE.




Приведем выражения для $M=G_2/T^2$ (т.е. для $j=8$).
{\footnotesize
\begin{multline*}{}
f_{14}(x)
= 2302911\,x^{14} - 10589454\,x^{13} + 25929072\,x
^{12} - 40992642\,x^{11} + 40560939\,x^{10} \\
\mbox{} - 19334556\,x^{9} - 8129184\,x^{8} + 21916564\,x^{7} -
17302275\,x^{6} + 5385490\,x^{5} \\
\mbox{} + 2042208\,x^{4} - 2822162\,x^{3} + 1247577\,x^{2} -
244584\,x + 32400
\end{multline*}
}%
Выражения $y=F_1(x)$ и $w=F_2(x)$ ($j=8$).
Представления многочленами степени $\le 13$ от $x$.
{\footnotesize
{\let\alpha x \begin{align*}{}
\mathit{y} :=
- {\displaystyle \frac {
175010896509363118040086341}{5631737006960999042017600}} \,\alpha
 ^{13}
&
+ {\displaystyle \frac {760112877954361925067757749}{
5631737006960999042017600}} \,\alpha ^{12}
\\
\mbox{} - {\displaystyle \frac {1750144201826669628973688007}{
5631737006960999042017600}} \,\alpha ^{11}  
&
\mbox{} + {\displaystyle \frac {2511833210232128226813555927}{
5631737006960999042017600}} \,\alpha ^{10}
\\
\mbox{} - {\displaystyle
\frac {1003878467681567742016705967}{2815868503480499521008800}}
\,\alpha ^{9} 
&
\mbox{}
+ {\displaystyle \frac {240879181713958481246898829}{
8447605510441498563026400}} \,\alpha ^{8}
\\
\mbox{}+ {\displaystyle
\frac {7403781632071474402565699143}{25342816531324495689079200}
} \,\alpha ^{7} 
&
\mbox{} - {\displaystyle \frac {27594158203874223011049838409}{
76028449593973487067237600}} \,\alpha ^{6}
\\
\mbox{}+ {\displaystyle
\frac {410494638054286810285781003}{2027425322505959655126336}}
\,\alpha ^{5} 
&
\mbox{} - {\displaystyle \frac {330826078483182090005453261}{
30411379837589394826895040}} \,\alpha ^{4}
\\
\mbox{} - {\displaystyle
\frac {3445072458977611467519796207}{50685633062648991378158400}
} \,\alpha ^{3}  
&
\mbox{} + {\displaystyle \frac {7587084783582358318419111269}{
152056899187946974134475200}} \,\alpha ^{2}
\\
\mbox{} - {\displaystyle
\frac {32695109346932045297414717}{2111901377610374640756600}} \,\alpha
&
\mbox{} + {\displaystyle \frac {185143322469969876688236}{
87995890733765610031525}}
\end{align*} }

{\let\alpha x \begin{align*}{}
\mathit{w} :=  - {\displaystyle \frac {
1576346706064295403636224151}{31537727238981594635298560}} \,
\alpha ^{13}
&
\mbox{}+ {\displaystyle \frac {2988198723955016496597743277
}{15768863619490797317649280}} \,\alpha ^{12} \\
\mbox{} - {\displaystyle \frac {13285340740238310997510265217}{
31537727238981594635298560}} \,\alpha ^{11} 
&
\mbox{} + {\displaystyle \frac {9257614313874971696870033121}{
15768863619490797317649280}} \,\alpha ^{10} \\
\mbox{} - {\displaystyle \frac {22633930581543413268015570611}{
47306590858472391952947840}} \,\alpha ^{9} 
&
\mbox{} + {\displaystyle \frac {8149391014299901771067458001}{
70959886287708587929421760}} \,\alpha ^{8} \\
\mbox{} + {\displaystyle \frac {99118813229806384154825380279}{
425759317726251527576530560}} \,\alpha ^{7} 
&
\mbox{} - {\displaystyle \frac {203481616765900685799134222281}{
638638976589377291364795840}} \,\alpha ^{6} \\
\mbox{} + {\displaystyle \frac {31059759863976878968413569399}{
170303727090500611030612224}} \,\alpha ^{5} 
&
\mbox{} - {\displaystyle \frac {4071162105662468344325556835}{
255455590635750916545918336}} \,\alpha ^{4} \\
\mbox{} - {\displaystyle \frac {39868741658792538217841910131}{
851518635452503055153061120}} \,\alpha ^{3} 
&
\mbox{} + {\displaystyle \frac {43885565296005598812277067281}{
1277277953178754582729591680}} \,\alpha ^{2}
\\
- {\displaystyle
\frac {474728963794058513671615399}{53219914715781440947066320}}
\,\alpha  
&
\mbox{} + {\displaystyle \frac {4723023511792833377635441}{
2956661928654524497059240}}
\end{align*} }
}


Представления в виде дробей.

{\footnotesize
{\let\alpha x \begin{multline*}{}
y=
 - 3((x - 1)(26244\,x^{10} - 107163\,x^{9} + 247860\,x^{8} -
349983\,x^{7} + 294081\,x^{6} - 116321\,x^{5} \\
\mbox{} - 38995\,x^{4} + 81535\,x^{3} - 52221\,x^{2} + 16668\,x
 - 2025))/( - 2700 - 1368\,x \\
\mbox{} - 25878\,x^{3} + 27987\,x^{2} + 19683\,x^{10} - 223074\,x
^{9} + 481140\,x^{8} - 327186\,x^{7} \\
\mbox{} - 114246\,x^{6} + 289122\,x^{5} - 123224\,x^{4})
\end{multline*} }%
}{\footnotesize
{\let\alpha x \begin{multline*}{}
w=
 - {\displaystyle \frac {1}{8}} (137781\,x^{11} - 216513\,x^{10}
 + 355752\,x^{9} - 943434\,x^{8} + 1340712\,x^{7} - 1027572\,x^{6}
\\
\mbox{} + 374878\,x^{5} + 60690\,x^{4} - 212989\,x^{3} + 209829\,
x^{2} - 96102\,x + 16200)/( - 2700 \\
\mbox{} - 1368\,x - 25878\,x^{3} + 27987\,x^{2} + 19683\,x^{10}
 - 223074\,x^{9} + 481140\,x^{8} \\
\mbox{} - 327186\,x^{7} - 114246\,x^{6} + 289122\,x^{5} - 123224
\,x^{4})
\end{multline*} }
}

\begin{REM}{} Введем в пространстве гомотетических классов инвариантных метрик
в $M/\sigma $ координаты $x=t_{b'}/t_b$, $y=t_{a'}/t_a$, $w=t_{b'}/t_a$
(см. Рис.3).
%
%
Система 
\cite[Eq. (1.1)]{2007}
после исключения
конформных множителей становится системой уравнений относительно $x,y,w$.
Программа MAPLE решила ее за секунду с четвертью
(понадобилось еще время на обработку решений).
\end{REM}

\section*{Дополнение 3. Положительно определенные эйнштейновы метрики на $M/\sigma $}

Пусть снова $M$ --- любое из флаговых пространств \eqref{eq:15}.


\begin{PROP}{} (a) На пространстве $M/ \sigma $ существуют положительно
оп\-ре\-де\-ленные $G$-инвариантные метрики Эйнштейна.
Число этих метрик с точ\-но\-стью до гомотетии
равно $2$ при $M=G_2/T^2$ (т.е. при $j=8$), и $6$ в остальных случаях.
\\
(b) Положительно определенные (соответственно, вещественные индефинитные)
инвариантные метрики Эйнштейна в $M/\sigma $,
отличные от метрик \eqref{eq:19},
отвечают положительным (соответственно, отрицательным)
корням неприводимого многочлена $f_{14}(x) \in \QQ[x]$.
\\
(c) Отрицательные корни $f_{14}(x)$ существуют только при $j\in \{36,60 \}$
(их по два).
\end{PROP}

\begin{proof}[Схема доказательства] То, что вещественные метрики
из (b) отвечают вещественным корням $f_{14}(x)$, следует
из леммы~\ref{LEM:D2} и предложения~\ref{PROP:1.4}.
Далее, можно проверить, что минимальными многочленами $\phi _i(y)\in \QQ[y]$ элементов
$y_i = F_i(x) \in \QQ[x]/(f_{14}(x))$, $i=1,2$, являются
некоторые
%
многочлены $14$ степени со знакочередующимися коэффициентами
\footnote{
Ю.Сакане \cite{Sa} аналогичным способом
проверил положительную опре\-де\-ленность
метрики Эйнштейна (iv) в $SU_4/T^3$
(см. выше, \S~\ref{sect:prim-pri-korr}, (iv)).
}.
Например, при $j=8:$
{\footnotesize
\begin{align*}{}
\phi _1(y)  &
=8290479600\,y^{14} - 52700208984\,y^{13} +
538432220079\,y^{12} - 1618143182406\,y^{11}
\\ &
\mbox{} + 2697711175278\,y^{10} - 3407189361930\,y^{9}
+ 3519687897648\,y^{8} - 2757517036290\,y^{7}
\\ &
\mbox{}  + 1539019875362\,y^{6} - 602718526838\,y^{5}
+ 166447276569\,y^{4} - 32207378760\,y^{3}
\\ &
\mbox{}  + 4184043012\,y^{2} - 327367656\,y  + 11573604  ,
\\
\phi _2(w) &
=1658095920\,w^{14} - 10281139488\,w^{13} + 29764527624 \,w^{12} - 53875631088\,w^{11}
\\ &
\mbox{} + 68589182916\,w^{10} - 65344826088\,w^{9} + 48093022782 \,w^{8} - 27694511532\,w^{7}
\\ &
\mbox{} + 12514329552\,w^{6} - 4426523804\,w^{5} + 1212637473\,w^{4} - 251395110\,w^{3}
\\ &
\mbox{} + 37706553\,w^{2} - 3645000\,w + 170586   .
\end{align*}{}
}
Поэтому $\phi _1(y)$ и $\phi _2(y)$ не имеют отрицательных корней,
что доказывает (b).
В следующем списке приводятся приближенные значения
всех вещественных корней $x_i$ многочлена $f_{14}(x)$ при различных $j:$
\begin{enumerate}\footnotesize
\item
$j=8$,  $x_1 = 1.0234 $, $x_2 = 1.0347  $,
\item
$j=18$, $x_1 = .3927 $, $x_2 = .4506   $,  $x_3 = .8181  $, $x_4 = 2.562  $,
\item
$j=24$, $x_1 =.3547  $, $x_2 = .4174  $,  $x_3 = .7419  $, $x_4 =  3.822  $,

\item
$j=36$, $x_1 = .3407   $, $x_2 =  .3843  $,  $x_3 =   .6845 $, $x_4 =  6.320   $,
$x_5 = -1.342   $,  $x_6 =    -1.285 $,
\item
$j=60$, $x_1 =   .3356 $, $x_2 =  .3618   $,  $x_3 =  .6472 $, $x_4 =   11.31  $,
$x_5 = -1.640   $,  $x_6 =   -1.547  $.
\end{enumerate}
Отсюда явствуют (a) и (c).
\end{proof}



\tiny

\end{document}